\documentclass[12pt]{article}
\usepackage{amsthm,amsfonts,amssymb,amscd}

\begin{document}

\begin{center}
\Large \bf Birationally rigid Fano hypersurfaces
\end{center}
\vspace{1cm}

\centerline{\large  A.V.Pukhlikov} \vspace{1cm}

\centerline{January 17, 2002}\vspace{1cm}

\parshape=1
3cm 10cm \noindent {\small \quad\quad\quad \quad\quad\quad\quad
\quad\quad\quad {\bf }\newline We prove that a smooth Fano
hypersurface $V=V_M\subset{\mathbb P}^M$, $M\geq 6$, is birationally
superrigid. In particular, it cannot be fibered into uniruled
varieties by a non-trivial rational map and each birational map
onto a minimal Fano variety of the same dimension is a biregular
isomorphism. The proof is based on the method of maximal
singularities combined with the connectedness principle of
Shokurov and Koll\' ar.} \vspace{1cm}

\centerline{}

\noindent 0. Introduction

0.1. The main result

0.2. Birationally rigid varieties

0.3. Log-canonical singularities

0.4. Maximal singularities

0.5. Scheme of the proof

0.6. Historical remarks and acknowledgements

\noindent 1. Projections and multiplicities

1.1. Singularities of subvarieties on hypersurfaces

1.2. How multiplicities change when we project

1.3. Singularities of the divisor $F$

\noindent 2. Log-canonical singularities

2.1. The connectedness principle of Shokurov and Koll\' ar

2.2. Log-canonical singularities of divisors in ${\mathbb P}^k$

\noindent 3. The main construction

3.1. The direct image of a maximal singularity

3.2. Comparing multiplicities

3.3. Computation of $\nu_E(F)$

\newpage

\section*{Introduction}

\subsection{The main result}

Fix an integer $M\geq 6$. Let ${\mathbb P}={\mathbb P}^M$ be the
complex projective space, $V=V_M\subset{\mathbb P}$ a smooth hypersurface
of degree $M$. Obviously, $V\subset{\mathbb P}$ is a Fano variety. By
the Lefschetz theorem $\mathop{\rm Pic}V={\mathbb Z} K_V={\mathbb Z} H$,
where $K_V=-H$, $H$ is the class of a hyperplane section on $V$.
The main result of the present paper is the following

{\bf Theorem.} {\it Variety $V$ is birationally superrigid.}

{\bf Corollary.} {\it
{\rm (i)} $V$ can not be fibered into uniruled varieties by a
non-trivial rational map.

{\rm (ii)} If $\chi\colon V -\, -\, \to V'$ is a birational map
onto a Fano variety $V'$  with ${\Bbb Q}$-factorial terminal
singularities such that $\mathop{\rm Pic} V'\otimes{\Bbb Q}={\Bbb
Q}K_{V'}$, then $\chi$ is a (biregular) isomorphism. In
particular, the groups of birational and
biregular self-maps coincide:
$$
\mathop{\rm Bir}V=\mathop{\rm Aut}V.
$$

{\rm (iii)} $V$ is non-rational.
}

Thus our main result is a higher-dimensional generalization of the
famous paper of Iskovskikh and Manin [IM], which made the starting point
of the modern birational geometry. Iskovskikh and Manin proved birational
superrigidity of three-dimensional quartics (the case $M=4$ in our
notations). The case $M=5$ was treated in [P1]. Superrigidity of
general (in the sense of Zariski topology) hypersurfaces $V$ was proved in
[P3]. Putting together [IM], [P1] and the present paper, we get
birational superrigidity of all smooth Fano hypersurfaces
$V=V_M\subset{\mathbb P}$ of dimension three and higher.

\subsection{Birationally rigid varieties}

Throughout this paper the ground field is assumed to be the field
${\mathbb C}$ of complex numbers. Recall (see [P2-P8]) that a smooth
Fano variety $X$ of
dimension $\geq 3$ with $\mathop{\rm rk}\mathop{\rm Pic} X=1$ is said
to be {\it birationally superrigid}, if for each birational map
$$
\chi\colon X-\,-\,\to X'
$$
onto a variety $X'$ of the same dimension, smooth in codimension
one, and each linear system $\Sigma'$ on $X'$, free in
codimension 1 (that is, $\mathop{\rm codim}\nolimits \mathop{\rm
Bs}\Sigma'\geq 2$), the inequality
\begin{equation}
\label{001} c(\Sigma,X)\leq c(\Sigma',X')
\end{equation}
holds, where $\Sigma=(\chi^{-1})_*\Sigma'$ is the strict transform
of $\Sigma'$ on $X$ with respect to $\chi$, and
$c(\Sigma,X)=c(D,X)$ stands for the {\it threshold of canonical
adjunction}
$$c(D,X)=\sup\{b/a|b,a\in{\mathbb Z}_+\setminus \{0\},
|aD+bK_X|\neq\emptyset\}$$ $D\in\Sigma$, and similarly for
$\Sigma'$, $X'$.

If for each triple $(X',\chi,\Sigma')$ there exists a birational
self-map $\chi^*\in\mathop{\rm Bir}X$ such that instead of the
inequality (\ref{001}) we have the estimate
\begin{equation}
\label{002} c(\Sigma^*,X)\leq c(\Sigma',X'),
\end{equation}
where the linear system
$$
\Sigma^*=(\chi\circ\chi^*)^{-1}_*\Sigma'=(\chi^*)^{-1}_*\Sigma
$$
is the strict transform of the linear system $\Sigma'$ with respect
to the composition
$$
\chi\circ\chi^*\colon X\stackrel{\chi^*}{-\,-\,\to}X\stackrel{\chi}{-\,-\,\to}X',
$$
then $X$ is said to be {\it birationally rigid}.

It is well known that a birationally
rigid variety $X$ can not be fibered into uniruled varieties of
a smaller dimension by a non-trivial rational map (this fact is an
immediate implication of the estimate (\ref{002}) combined with
the assumption that $\mathop{\rm Pic} X\cong{\mathbb Z}$, see [P3]),
and that if $\chi\colon X -\, -\, \to X'$ is a birational map
onto a Fano variety $X'$  with ${\Bbb Q}$-factorial terminal
singularities such that $\mathop{\rm Pic} X'\otimes{\Bbb Q}={\Bbb
Q}K_{X'}$, then $X'$ is (biregularly) isomorphic to $X$, see, for
instance, [P7]. In particular, birationally rigid varieties are
non-rational.

Birational rigidity and superrigidity can also be defined for fibrations
into Fano varieties, see [P4,P6] and also [S]. However, this
relative version of rigidity will not be discussed in the present
paper.

Birational (super)rigidity is one of the most striking and
mysterious phenomena of higher-dimensional algebraic geometry.
Informally speaking, rigidity means that algebraic varieties with
no global non-zero differential forms (rationally connected
varieties) behave as if there were plenty of differential
forms on them. Indeed, the main result of the present paper
implies that for any smooth hypersurface $V'\subset{\mathbb P}$
any birational map
$$
\chi\colon V-\,-\,\to V'
$$
is a (projective) isomorphism. If the degree of the hypersurface
$V$ were higher than $M+1$, this claim would have been obvious:
just compare the very ample canonical linear system
$|K_V|$ and the system $|K_{V'}|$. However, in our case the canonical
linear system $|K_V|$ (together with all the pluricanonical linear systems
$|mK_V|$, $m\geq 1$) is empty, so that we have nothing to compare.
And nevertheless the hypersurface $V$ behaves more like a variety
of general type than the projective space, its close neighbour in
the class of Fano varieties. What are the deep reasons underlying
the phenomenon of birational (super)rigidity?  An answer to this
question is yet to be found.

\subsection{Log-canonical singularities}

Recall certain definitions of the theory of log-minimal models.
For simplicity assume $X$ to be a smooth algebraic variety,
$D=\sum\limits_{i\in I_D}d_iD_i$ a ${\Bbb Q}$-divisor, where
$0<d_i<1$, $D_i\subset X$ are pair-wise distinct irreducible
hypersurfaces.

{\bf Definition 1.} An irreducible subvariety $W\subset X$ is
called a {\it log-center} of the pair $(X,D)$, if for some
resolution of singularities $f\colon Y\to X$ of the pair $(X,D)$
with the set $\{ E_i\subset Y\,|\, i\in I_f\}$ of exceptional
divisors we get
$$
K_Y+\widetilde D= f^*(K_X+D)+\sum_{i\in I_f}e_iE_i,
$$
where for some $i\in I_f$
$$
f(E_i)=W\quad\mbox{and}\quad e_i\leq -1.
$$
The set of log-centres is not an invariant of the pair $(X,D)$.
However, the set-theoretic union
$$
\mathop{\rm LC} (X,D)=\mathop{\bigcup}\limits_{e_i\leq -1}f(E_i)
$$
of all log-centres depends on the pair $(X,D)$ only. For an
arbitrary point $x\in X$ denote by the symbol
$$
\mathop{\rm LC} (X,D,x)
$$
the connected component of the closed algebraic set $\mathop{\rm
LC} (X,D)$, containing $x$.

{\bf Definition 2.} The pair $(X,D)$ is {\it log-terminal at the
point $x$ in dimension} $l\geq 1$, if
$$
\dim \mathop{\rm LC} (X,D,x)\leq l-1.
$$
The pair $(X,D)$ is {\it log-terminal at the point} $x$, if $x\not\in
\mathop{\rm LC} (X,D)$, or equivalently $\mathop{\rm LC}
(X,D,x)=\emptyset$.

{\bf Definition 3.} The pair $(X,D)$ is {\it log-canonical at the
point $x$ in dimension} $l\geq 1$, if the pair
$$
(X,\frac{1}{1+\varepsilon} D)
$$
is log-terminal at the point $x$ in dimension $l$ for all
sufficiently small $\varepsilon>0$.

In a similar way we define the property of being log-canonical at
the point $x$.

Now consider the pair $({\mathbb P}^k,D)$. We define the {\it
degree} of the ${\Bbb Q}$-divisor $D$ in a natural way as the
number
$$
\mathop{\rm deg} D=\sum_{i\in I_D}d_i \mathop{\rm deg} D_i,
$$
so that $D\equiv (\mathop{\rm deg} D)H$ in $A^1{\Bbb
P}^k\otimes{\Bbb Q}$.

{\bf  Proposition 1.} {\it If the pair $({\mathbb P}^k,D)$ is
log-terminal (respectively, log-canonical) at the point
$x\in{\mathbb P}^k$ in dimension $l\geq 1$, where the integer $l$
satisfies the inequality
\begin{equation}
\label{1a}
l+\mathop{\rm deg} D\leq k+1,
\end{equation}
then it is log-terminal (respectively, log-canonical) at the
point $x$.}

{\bf Proof} is given in Section 2.

\subsection{Maximal singularities}

Assume that there exists a birational map
$$
\chi\colon V-\,-\,\to V'
$$
onto a variety $V'$ of the same dimension, smooth in codimension
one, and a moving linear system $\Sigma'$ on $V'$ such that for
the strict transform $\Sigma$ of the system $\Sigma'$ the
inequality
$$
c(\Sigma,V)> c(\Sigma',V')
$$
holds. Obviously, $c(\Sigma,V)=n\geq 1$, where
$\Sigma\subset |nH|$.

{\bf Proposition 2.} {\it In the assumtions above there exists a
geometric discrete valuation $\nu$ (that is, a valuation which is
divisorial on a certain model $V^{\sharp}$ of the variety $V$)
satisfying the Noether-Fano inequality
\begin{equation}
\label{0051}
\nu (\Sigma)> n\cdot \mathop{\rm discrepancy}(\nu).
\end{equation}
}

{\bf Definition 4.} A geometric discrete valuation
satisfying the Noether-Fano inequality (\ref{0051}) is called a
{\it maximal singularity} of the linear system $\Sigma$.

More explicitly, (\ref{0051}) means that there is a birational
morphism
$$
\varphi\colon V^{\sharp}\longrightarrow V
$$
of projective varieties and a prime divisor $E\subset V^{\sharp}$
such that $E\not\subset\mathop{\rm Sing}V^{\sharp}$ and
$$
\mathop{\rm ord}\nolimits_E(\varphi^*\Sigma)>na(E),
$$
where $a(E)$ stands for the discrepancy of $E$ with respect to
the model $V$.

{\bf Proof of Proposition 2}  can be found in many papers; see,
for instance, [C,CPR,I,IM,IP,P2-P6]. It is quite elementary.

It is easy to see that for the
centre $B=\varphi(E)$ of the maximal singularity $\nu$
on $V$ we have the inequality
\begin{equation}
\label{0053}
\mathop{\rm mult}\nolimits_B\Sigma>n.
\end{equation}

{\bf Remark.} If the birational morphism $\varphi$ above is just
the blow up of an irreducible subvariety $B\subset V$, then the
condition of the exceptional divisor $E\subset V^{\sharp}$ being a
maximal singularity takes an especially simple form:
$$
\mathop{\rm mult}\nolimits_B\Sigma>n(\mathop{\rm codim}B-1).
$$
This particular case of the Noether-Fano inequality is called the
{\it Fano inequality}, and the subvariety $B\subset V$ in this case
is called a {\it maximal subvariety} (point, curve, surface,\dots)
of the linear system $\Sigma$. Alas, the general theory makes sure
that a maximal singularity does exist but says nothing about
existence of its nicest particular case, a maximal subvariety. And
indeed, there are examples (see, for instance, the recent preprints
[G1-G3] of Gizatullin) when there is no maximal subvariety whereas
one can easily find a maximal singularity making {\it a few} blow
ups. If a maximal singularity is obtained by means of a few (that
is, more than one) blow ups, it is said  to be
{\it infinitely near}.

In other words, what the general theory asserts is existence either
of a maximal subvariety or of an infinitely near maximal singularity.
Quite predictably, it is the infinitely near case that makes the
work really hard and requires sometimes considerable effort to cope
with. Anyway, the reader should keep this alternative in mind to
avoid misunderstanding which might be generated by [G1-G3].

\subsection{Scheme of the proof}

Assume that the variety $V$ is not birationally superrigid. By the
previous section and the definition of superrigidity this implies
that there exists a moving linear system $\Sigma\subset |nH|$
with a maximal singularity $\nu$. Now, the following fact is true:

{\bf Proposition 3.} {\it For any curve $C\subset V$ the
estimate
$$
\mathop{\rm mult}\nolimits_C \Sigma\leq n
$$
holds.}

{\bf Proof:} see [P2,P3].

Therefore by (\ref{0053}) the centre of the discrete valuation
$\nu$ on $V$ is a
point $x$. Let $D_1,D_2\in\Sigma$ be general divisors,
$$
Z=(D_1\circ D_2)
$$
the effective algebraic cycle of the scheme-theoretic
intersection of the divisors $D_1$ and $D_2$.

Take a general line $L\subset{\mathbb P}$ and let
$$
\pi\colon{\mathbb P}-\,-\,\to {\mathbb P}^{M-2}
$$
be the corresponding linear projection. We may assume that
$\mathop{\rm Supp}Z\cap L=\emptyset$,
so that $\pi_L|_V$ is regular in a neighborhood of
the set $\mathop{\rm Supp} Z$. Consider the direct image
$$
F=(\pi_L)_* Z.
$$
It is an effective divisor on ${\mathbb P}^{M-2}$. By the Lefschetz
theorem the multiplicity of each component of the cycle $Z$ is
not higher than $n^2$, since $Z\equiv n^2H^2$. Consequently, we
may assume that this is true for $F$, either. Obviously,
$\mathop{\rm deg} F=Mn^2$.

The following fact is crucial in our arguments. Set $y=\pi_L(x)$
to be the image of the point $x$ in ${\Bbb P}^{M-2}$.

{\bf Proposition 4.} {\it The pair
$({\Bbb P}^{M-2},\frac{1}{2n^2}F)$ is {\rm (i)} log-canonical
at the point $y$ in dimension 2, but {\rm (ii)} not log-canonical
at the point $y$.}

{\bf Proof} of part (i) is given in Section 1, of part (ii) in
Section 3.

{\bf End of the proof of the theorem.} By Proposition 1 we get the
estimate
$$
2+\mathop{\rm deg} (\frac{1}{2n^2}F)=2+\frac{M}{2}>M-1,
$$
whence $M<6$: a contradiction. Hence our assumption that $V$ is
not birationally superrigid is false. Q.E.D.

\subsection{Historical remarks and acknowledgements}

As we have mentioned above, birational superrigidity of a general
hypersurface $V=V_M\subset{\mathbb P}^M$ was proved in [P3]. Later
superrigidity was proved for general hypersurfaces with isolated
non-degenerate singularities [P8]. In 1998 at the Satellite
conference on algebraic geometry in Essen Corti suggested to use
the Shokurov-Koll\' ar connectedness principle [Sh,K] for proving
birational rigidity. However, the first version of the paper [C]
contained a mistake coming from an overenthusiastic
generalization of the connectedness principle. Using the wrong
argument of Corti, Cheltsov in a few months produced a ``proof''
of birational superrigidity of an arbitrary smooth hypersurface
$V=V_M\subset{\mathbb P}^M$. However, the estimates that he used in
his arguments were too strong to be realistic. In December 1998
simple counterexamples to these estimates were constructed (by the
author of the present paper and Grinenko). As a by-product, the
mistake was discovered in the first version of Corti's paper [C].
Unfortunately, after this mistake had been corrected, it turned
out that the new technique made it possible mainly just to give
new proofs of the theorems already proved by the classical
methods. In the papers [Ch1,2] Cheltsov proved superrigidity of
an arbitrary smooth hypersurface for $M\in\{5,6,7,8\}$; however,
some of his arguments are rather doubtful (see a detailed
exposition in [I]).

At the same time, the papers of Park and Cheltsov [Pk,ChPk,Ch3]
left some hope that new stronger estimates for multiplicities of
divisors on hypersurfaces could be obtained by means of the
Shokurov-Koll\' ar connectedness principle. Indirectly this hope
was supported by one argument of Corti's [C] (considerably
simplified by Cheltsov [Ch3]) that excluded a maximal singularity
over a three-fold non-degenerate double point. These ideas were
further developed in [CM].

In this paper we use both the classical method of maximal
singularities and the connectedness principle. A combination of
various ideas makes it possible to prove, at long last, the
superrigidity of Fano hypersurfaces
without the annoying assumption of general position
(regularity, see [P3]). However, for other principal classes of
higher-dimensional Fano varieties the methods developed here are
still not strong enough.

The author was supported by Alexander von Humboldt Foundation, Russian
Foundation of Basic Research and INTAS and by Science
Support Foundation grant for young researches.

\newpage

\section{Projections and multiplicities}

The aim of this section is to prove part (i) of Proposition 4.
The arguments break into two components: firstly, we prove an
estimate for the multiplicity of the cycle $Z$ along an arbitrary
surface, secondly, we extend these estimates to the divisor $F$ in
${\mathbb P}^{M-2}$.

\subsection{Singularities of subvarieties on hypersurfaces}

To begin with, let us consider the following general situation.
Let $X\subset{\mathbb P}$ be an arbitrary hypersurface, $D$ an
effective divisor that is cut out on  $X$ by a hypersurface
of degree $n$ in ${\mathbb P}$, that is, $D\in|nH|$. It was proved
in [P2,P3] that for any curve $C\subset
X\setminus\mathop{\rm Sing} X$ the following estimate holds:
\begin{equation}
\label{007}
\mathop{\rm mult}\nolimits_CD\leq n.
\end{equation}
An immediate implication of this claim is Proposition 3 above.
However, the estimate (\ref{007}) is a particular case of the
following general fact which is crucial for our proof of
Proposition 4, (i).

{\bf Proposition 5.} {\it Let $W$ be an effective cycle of
codimension $k$, $W\sim mH^k$ in the group of cycles
$A_{M-k-1}X$, $S\subset W$ an irreducible subvariety of dimension
$k$, where $2k+1<M$. Assume that $S\cap\mathop{\rm
Sing}X=\emptyset$. Then the following estimate is true:
$$
\mathop{\rm mult}\nolimits_S W\leq m.
$$
}

{\bf Proof.}  We use the same construction as in [P2]: for an irreducible
subvariety $T\subset X$ of dimension $k$,
$T\cap\mathop{\rm Sing}X=\emptyset$, and a
point $x\in{\mathbb P}\setminus X$ set $C(x,T)\subset{\mathbb P}$ to be
the cone with the vertex at the point $x$ and the base $T$.
Obviously,
\begin{equation}
\label{008}
X\cap C(x,T)=T\cup R(x,T),
\end{equation}
where $R(x,T)$ is the residual closed set of the same dimension as $T$. It is
easy to see that for a sufficiently general point $x$ the closed set
$R(x,T)$ is an irreducible subvariety and the presentation (\ref{008}) holds
not only in the set-theoretic but also in the scheme-theoretic sense.

Now let $B\subset X$ be any closed algebraic set.

{\bf Lemma 1.} {\it Assume that the point $x\in{\mathbb P}\setminus X$ is
sufficiently general. Then: {\rm (i)} if $\mathop{\rm codim}\nolimits_X B>
\mathop{\rm dim} T$, then
$$
R(x,T)\cap B\subset T\cap B;
$$
{\rm (ii)} if $\mathop{\rm codim}\nolimits_X B=\mathop{\rm dim} T$, then the
intersection $R(x,T)\cap B$ is zero-dimensional outside $T$.}

An (easy) proof of the lemma is given below.

Let us prove Proposition 5. Take a general $k$-uple of points
$$
(x_1,\dots,x_k)\in ({\mathbb P}\setminus X)^k
$$
and construct by induction on $i=1,\dots,k$ the following sequence of irreducible
$k$-dimensional subvarieties, starting with $R_0=S$:
$$
R_i=R_i(x_1,\dots,x_i)=R(x_i,R_{i-1}).
$$

{\bf Lemma 2.} {\it Assume that the $k$-uple $(x_1,\dots,x_k)$ is sufficiently
general. Then the intersection
$$
R_k\cap S
$$
consists of precisely $\mathop{\rm deg}R_k=(\mathop{\rm deg}X-1)^k
\mathop{\rm deg}S$ distinct points.}

Putting off proof of the lemma for a while, let us complete the proof
of Proposition 5. By Lemma 1 the intersection $R_k\cap W$ is zero-dimensional
outside $S$. Therefore
\begin{equation}
\label{0014}
\begin{array}{c}\displaystyle
m\mathop{\rm deg}R_k=(R_k\cdot W)\geq \\  \\
\displaystyle
\geq \sum_{p\in R_k\cap S}(R_k\cdot W)_p\geq \\ \\
\displaystyle
\geq \mathop{\rm mult}\nolimits_S W\cdot \mathop{\rm deg}R_k.
\end{array}
\end{equation}

This estimate implies immediately Proposition 5. Q.E.D.

{\bf Proof of Lemma 1.} For two closed sets $A,B\subset{\mathbb P}$ define
$$
J(A,B)=\overline{\mathop{\bigcup}\limits_{A\ni p\neq q\in
B}L_{p,q}}
$$
to be the closure of the set, which is swept out by all lines
connecting points of $A$ and $B$. Obviously,
\begin{equation}
\label{0017}
\dim J(A,B)\leq \dim A+\dim B+1.
\end{equation}
Now if $\mathop{\rm codim}\nolimits_X B>
\mathop{\rm dim} T$ then the set $J(T,B)$ has positive codimension
in ${\mathbb P}$. If a point
$$
y\in R(x,T)\cap B
$$
does not belong to $T$, then obviously $x\in J(T,B)$. This proves
part (i). If $\mathop{\rm codim}\nolimits_X B=\mathop{\rm dim} T$
then, generally speaking, we have the equality $J(T,B)={\mathbb P}$.
In this case denote by ${\cal C}\subset{\mathbb P}\times{\mathbb P}
\times{\mathbb P}$ the closure of the set of all collinear triples
$(u,v,w)$ with the natural projections $\pi_1,\pi_2,\pi_3$ onto
the factors ${\mathbb P}$. Now
$$
\pi_3\colon {\cal C}_{T,B}={\cal C}\cap \pi^{-1}_1(T)\cap \pi^{-1}_2(B)
\to J(T,B)={\mathbb P}
$$
is a surjective morphism. Since $\mathop{\rm dim}{\cal C}_{T,B}=
\mathop{\rm dim}{\mathbb P}$, it is generically finite. Thus for a
general point $x\in {\mathbb P}\setminus X$ the inverse image
$\left(\pi_3|_{{\cal C}_{T,B}}\right)^{-1}(x)$ is a finite set of points.
Now as above, if a point $y\in R(x,T)\cap B$ does not belong to $T$, then
there exists a point $z\in T$ such that the triple $(z,y,x)$ is collinear.
In other words,
$$
(z,y,x)\in\left(\pi_3\left|_{{\cal C}_{T,B}}\right.\right)^{-1}(x),
$$
which is what we need. Q.E.D. for Lemma 1.

{\bf Proof of Lemma 2.} Let us consider an arbitrary
$k$-dimensional irreducible subvariety
$T\subset X$, lying outside the singular locus, that is,
$T\cap\mathop{\rm Sing}X=\emptyset$, and compute
the intersection $R(x,T)\cap T$ explicitly. Let
$$
\pi_x\colon X\to{\mathbb P}^{M-1}
$$
be the regular map induced by the linear projection ${\mathbb P}-\,-\,\to
{\mathbb P}^{M-1}$ from the point $x$. Denote by $D(x)\subset X$ the branch
divisor of the morphism $\pi_x$.

{\bf Lemma 3.} {\it If the point $x$ is sufficiently general, then
set-theoretically
$$
R(x,T)\cap T=D(x)\cap T.
$$}

{\bf Proof.} Recall that by assumption $T$ is contained in the smooth
part of $X$. At a point $y\in T\setminus D(x)$ the differential
$$
(d\pi_x)_y\colon T_yX\to T_{\pi_x(y)}{\mathbb P}^{M-1}
$$
is an isomorphism (so that $\pi_x$ is a local isomorphism in the complex
analytic sense). Taking into account that the map
$\pi_x\colon T\to\pi_x(T)$ is bijective (since $k<(M-1)/2$ by assumption),
we see that near the point $y$ the inverse image
$\pi^{-1}_x(\pi_x(T))$ coincides with $T$. This implies that
$y\not\in R(x,T)$, so that set-theoretically
$$
R(x,T)\cap T\subset D(x)\cap T.
$$
Let us show that this inclusion is actually an equality. Since the point
$x$ is general we may assume that $T\cap D^{\circ}(x)$ is dense in
$T\cap D(x)$, where $D^{\circ}(x)$ consists of the points where  $\pi_x$
is simply branched, $D^{\circ}(x)\subset D(x)$ is an open subset. Moreover,
for a suitable set $(z_1,\dots,z_{M-1})$ of local parameters at a point
$y\in T\cap D^{\circ}(x)$ the morphism $\pi_x$ looks as follows:
$$
(z_1,\dots,z_{M-1})\mapsto (z_1^2,z_2,\dots,z_{M-1})
$$
and we may assume that $T\neq \pi_x^{-1}(\pi_x(T))$ near $y$. In the
complex analytic setting, this means that $T$ is not $\sigma$-invariant
where
$$
(z_1,\dots,z_{M-1})\mapsto (-z_1,z_2,\dots,z_{M-1})
$$
is the local involution generated by the double cover. It is easy to see
that near the point $y$
$$
T\cap \sigma(T)=T\cap D(x),
$$
where $D(x)=\{z_1=0\}$. Since obviously
$$
\pi_x^{-1}(\pi_x(T))=T\cup \sigma(T)
$$
in a small complex analytic neighborhood of the point $y$, we get the
inclusion
$$
D^{\circ}(x)\cap T\subset R(x,T)\cap T.
$$
Taking the closure, we obtain our claim. Q.E.D. for Lemma 3.

Note that we used the complex analytic language in the proof just for
simplicity and clarity of exposition; the fact can be proved in the
obvious way in the purely algebraic setting.

Now take a system of homogeneous coordinates $(z_0\colon
z_1\colon \dots \colon z_M)$ on ${\mathbb P}$ and let $F(z_*)$ be an
equation of the hypersurface $X$. If $(x_0\colon x_1\colon \dots
\colon x_M)$ are the coordinates of the point $x$, then
$$
T\cap D(x)=\left\{ \left.\sum^M_{i=0}\frac{\partial F}{\partial
z_i}x_i\right|_T=0\right\}.
$$
Since $T\cap\mathop{\rm Sing}X=\emptyset$, the linear system
\begin{equation}
\label{1} \left|\left.\sum^M_{i=0}\lambda_i\frac{\partial
F}{\partial z_i}\right|_T=0\right|
\end{equation}
of divisors on $T$ is base point free. In particular, for a suffiently
general point $x$ the divisor $D(x)\cap T$ is
irreducible and reduced.

Here we come back to the geometric situation described in
Lemma 2. Recall that $k<\mathop{\rm dim}X/2$.

{\bf Lemma 4.} {\it Assume that the $k$-uple of points
$(x_1,\dots,x_k)$ is sufficiently general. Then for any
$0\leq i<j<l\leq k$ we have
\begin{equation}
\label{009}
R_i\cap R_l\subset R_j.
\end{equation}
}

{\bf Proof.} It is enough to consider the case $l=j+1$, but
here our claim follows from Lemma 1 (i), since $R_l=R(x_l,R_j)$
and $\mathop{\rm dim}R_i=\mathop{\rm dim}R_j=
k<\mathop{\rm dim} X/2$. Q.E.D.

Iterating (\ref{009}) $k$ times we obtain the following
set-theoretic equality
\begin{equation}
\label{0010}
R_k\cap S=\mathop{\bigcap}\limits^k_{i=0}R_i.
\end{equation}
Now since by Lemma 3 $R_{i+1}\cap R_i=D(x_{i+1})\cap R_i$, we
derive from (\ref{0010}) that set-theoretically
$$
\begin{array}{c}
\displaystyle
R_k\cap S= D(x_k)\cap \mathop{\bigcap}\limits^{k-1}_{i=0}R_i=\dots \\ \\
\displaystyle
\dots = \mathop{\bigcap}\limits^{k}_{i=j+1}D(x_i)\cap
\mathop{\bigcap}\limits^{j}_{i=0}R_i=\dots \\ \\
\displaystyle
\dots = \mathop{\bigcap}\limits^{k}_{i=1}D(x_i)\cap S.
\end{array}
$$
But the divisors $D(x_*)|_S$ are arbitrary divisors of a base point free
linear system. Therefore the intersection $R_k\cap S$ consists of
$$
(D(x)^k\cdot S)
$$
distinct points. It remains to add that $D(x)\in |(\mathop{\rm deg}X-1)H|$,
where $H$ is the class of a hyperplane section of the hypersurface
$X$. Q.E.D. for Lemma 2.

Proof of Proposition 5 is now complete.

{\bf Remark.} Let us explain why the number of points in the intersection
$R_k\cap S$ turns out to be precisely what we need to obtain our
estimate. Let us look at the cone $C=C(x,T)$ once again. Take a
resolution of singularities
$$
\tau\colon \widetilde T\to T
$$
and blow up the point $x$:
$$
\varphi\colon \widetilde{\mathbb P}\to{\mathbb P}.
$$
Denote by
$$
\pi\colon \widetilde {\mathbb P}\to {\mathbb P}^{M-1}
$$
the regularized projection from the point $x$ onto a general hyperplane.
The morphism $\pi$ presents $\widetilde {\mathbb P}$ as a locally trivial
${\mathbb P}^1$-bundle over ${\mathbb P}^{M-1}$. Taking the composition
$$
\beta\colon\widetilde T\stackrel{\tau}{\longrightarrow}T
\hookrightarrow \widetilde {\mathbb P} \stackrel{\pi}{\longrightarrow}
{\mathbb P}^{M-1}
$$
(where the inclusion $T\hookrightarrow\widetilde{\mathbb P}$ makes sense
since $x\not\in T$), look at the fiber product
$$
\widetilde C=\widetilde C(x,T)=
\widetilde T\mathop{\times}\nolimits_{{\mathbb P}^{M-1}}\widetilde {\mathbb P}.
$$
Write down the natural projection onto the first factor as
$p\colon \widetilde C\to \widetilde T$. It is a locally trivial
${\mathbb P}^1$-bundle over $\widetilde T$. Of course, the image of
$\widetilde C$ in $\widetilde {\mathbb P}$ is equal to the strict transform
of the cone $C$ on $\widetilde {\mathbb P}$. Thus we obtain $\widetilde C$
from $C$ by means of two operations: first, we blow up the vertex
$x$ of the cone, second, we pull back the resolution $\tau$ to the
corresponding locally trivial ${\mathbb P}^1$-bundle over $\pi(T)$. Denote
this map by
$$
\mu\colon \widetilde C\to C.
$$
Since $\widetilde C$ is a ${\mathbb P}^1$-bundle over $\widetilde T$, we
get
\begin{equation}
\label{0011}
\mathop{\rm Pic}\widetilde C=
{\mathbb Z}E\oplus p^*\mathop{\rm Pic}\widetilde T,
\end{equation}
where $E=\mu^{-1}(x)$ is the exceptional section,
$E\cong \widetilde T$. Let $T^{\sharp}=\mu^{-1}(T)$ be the
``base'' section (that is, coming from the original base
of the cone) and $H^{\sharp}=\mu^* H$ be the pull back of the
hyperplane section of the cone $C\subset{\mathbb P}$. Since both
$T^{\sharp}$ and $H^{\sharp}$ are sections of the bundle, by
(\ref{0011}) we get
$$
T^{\sharp}\sim H^{\sharp}+p^* Q,
$$
where $Q\in\mathop{\rm Pic}\widetilde T$. Intersecting with the
exceptional section, however, we get
\begin{equation}
\label{0012}
p^* Q|_E=0
\end{equation}
and thus $Q=0$, because the composition
$$
\mathop{\rm Pic}\widetilde T
\stackrel{p^*}{\longrightarrow}
\mathop{\rm Pic}\widetilde C
\stackrel{{\cdot}|_E}{\longrightarrow}
\mathop{\rm Pic} E
$$
is an isomorphism: it is induced by the natural isomorphism
$p\colon E\cong \widetilde T$. The crucial fact that implies
(\ref{0012}) is the obvious observation that both $T\subset C$
and $H\subset C$ (for a general hyperplane section) do not
contain the vertex of the cone and for this reason
$T^{\sharp}|_E=H^{\sharp}|_E=0$. We conclude that
\begin{equation}
\label{0013}
T^{\sharp}\sim H^{\sharp}.
\end{equation}
In other words, on the cone $C$ the divisor $T=\mu_* T^{\sharp}$
is rationally equivalent to a hyperplane section! In particular,
we deduce from (\ref{0013}) that
\vspace{1cm}

\parshape=1
3cm 10cm \noindent {\it
the intersection $R(x,T)\cap T$ is rationally equivalent
on $R(x,T)$ to a hyperplane section of $R(x,T)$.
} \vspace{1cm}

This is true for any subvariety $T\subset X$. Now let us count
the points of intersection of $R_k$ and $S$. By (\ref{0010}), it is
the same as to count the points in the set
$$
\mathop{\bigcap}\limits^{k}_{i=0}R_i=
R_0\cap R_1\cap R_2\cap\dots\cap R_k.
$$
Denote by $Q_j$, $j=0,1,\dots,k$, the intersection
$$
\mathop{\bigcap}\limits^{k}_{i=j}R_i.
$$
Since obviously $Q_j\subset R_j$, we obtain from (\ref{0013}) that
\vspace{1cm}

\parshape=1
3cm 10cm \noindent {\it
$Q_{j-1}=R_{j-1}\cap Q_j$ is rationally equivalent on $Q_j$ to a
hyperplane section of $Q_j$.
} \vspace{1cm}

Therefore $Q_0=R_k\cap S$ is rationally equivalent to the intersection
of $R_k$ with a plane $P\subset{\mathbb P}$ of codimension $k$,
that is, to a zero-dimensional scheme of degree
$\mathop{\rm deg}R_k$. To prove Proposition 5, we need only
this fact, not the precise value of $\mathop{\rm deg}R_k$.

A drawback of this argument is that a zero-dimensional scheme may
contain multiple points, so that to make this line of arguments work
one should check that the input of each point of the intersection
$R_k\cap S$ into the estimate (\ref{0014}) is correct (that is, not
less than its multiplicity in $R_k\cap S$). This is of course true
and not so hard to verify but nevertheless requires some extra rather
tedious considerations. So we prefer the straightforward computation
performed above leaving the present argument just as a useful
explanation.

\subsection{How multiplicities change when we project}

Abusing our notations, we write ${\mathbb P}$ instead of ${\Bbb
P}^k$. Let $Q\subset{\mathbb P}$ be an irreducible subvariety of
codimension $\mathop{\rm codim}Q\geq 2$, $q\in Q$ a point. For a
point $a\in{\mathbb P}$ of general position the morphism
$$
\pi_a\colon Q\to R\subset{\mathbb P}^{k-1},
$$
induced by the linear projection $\pi_a\colon{\mathbb P}
-\,-\,\to{\mathbb P}^{k-1}$ from the point $a$, is birational. Set
$r=\pi_a(q)$.

{\bf Proposition 6.} {\it The multiplicities coincide:
\begin{equation}
\label{2}
\mathop{\rm mult}\nolimits_rR=\mathop{\rm
mult}\nolimits_qQ
\end{equation}
if and only if the following two conditions are satisfied:

{\rm (i)} the line
$L_{a,q}=\overline{\pi^{-1}_a(\pi_a(q))}\subset{\mathbb P}$ that
connects the points $a$ and $q$ has no other points of
intersection with $Q$:
$$
L_{a,q}\cap Q=\{q\},
$$

{\rm (ii)} the line $L_{a,q}$ is not tangent to $Q$ at the point
$q$:
$$
L_{a,q}\not\subset T_qQ.
$$
}

{\bf Proof.} For a general plane $P\subset {\mathbb P}^{k-1}$, $P\ni
r$, of dimension $\mathop{\rm codim}\nolimits_{\mathbb P}Q-1$ we have
$$
\mathop{\rm mult}\nolimits_rR=(P\cdot R)_r.
$$
Let $P^+=\overline{\pi^{-1}_a (P)}$ be the inverse image of the
plane $P$ in ${\mathbb P}$. We get
$$
(P\cdot R)_r= \sum_{y\in\pi^{-1}_a(r)}(P^+\cdot Q)_y= (P^+\cdot
Q)_q+\sum_{\scriptstyle
\begin{array}{c}
\scriptstyle y\in\pi^{-1}_a(r),\\ \scriptstyle y\neq q
\end{array}
}(P^+\cdot Q)_y,
$$
where obviously $(P^+\cdot Q)_q\geq\mathop{\rm mult}\nolimits_q
Q$. Thus the equality (\ref{2}) holds if and only if the condition
(i) is satisfied and $(P^+\cdot Q)_q=\mathop{\rm mult}\nolimits_q
Q$ into the bargain. Let
$$
\varphi\colon X\to{\mathbb P}
$$
be the blow up of the point $q$, $E\cong {\mathbb P}^{k-1}$ the
exceptional divisor, $\widetilde Q$ and $\widetilde P^+$ the
strict transforms of $Q$ and $P^+$ on $X$. It is easy to see that
for a general plane $P$ the intersection
$$
\widetilde Q\cap\widetilde P^+
$$
is zero-dimensional, so that by the elementary
formulae of intersection theory [F] we get
$$
(P^+\cdot Q)_q=\mathop{\rm mult}\nolimits_q Q+ \sum_{y\in
\widetilde Q\cap \widetilde P^+ \cap E} (\widetilde Q\cdot
\widetilde P^+)_y.
$$
Finally, (\ref{2}) holds if and only if for a general plane $P$
the intersection
$$
\widetilde Q\cap \widetilde P^+\cap E
$$
is empty. But $\widetilde P^+\cap E$ is an arbitrary plane of
dimension $\dim P$, containing the point $\widetilde L_{a,q}\cap
E=T_q L_{a,q}$. Consequently, the latter condition takes the form
$$
\widetilde L_{a,q}\cap E\not\in \widetilde Q\cap E,
$$
that is, $L_{a,q}\not\subset T_qQ$, as we claimed. Q.E.D.

{\bf Corollary 1.} {\it If $\mathop{\rm codim}Q\geq 3$, then for
any (possibly reducible)
curve $C\subset Q$ for a sufficiently general point
$a\in{\mathbb P}$ the projection
$$
\pi_a\colon Q\to R
$$
is one-to-one and preserves multiplicity in a neighborhood of the
curve $C$. More precisely, there exists an open subset $U\subset
Q$, that contains entirely the curve $C$ (the complement
$Q\setminus U$ is of codimension 2 in $Q$), such that
$(\pi_a|_Q)^{-1}(\pi_a(U))=U$, the map $\pi_a\colon U\to
\pi_a(U)\subset R$ is bijective and for each point $z\in U$ the
following equality holds:}
$$
\mathop{\rm mult}\nolimits_{\pi_a(z)}R=\mathop{\rm
mult}\nolimits_zQ.
$$

{\bf Proof.} We use the notation $J(A,B)$ for any two closed
subsets $A,B\subset{\mathbb P}$, introduced in Sec. 1.1. Besides,
for a closed subset $A\subset Q$ let
$$
T(A,Q)=\overline{\mathop{\bigcup}\limits_{p\in A}T_pQ}
$$
be the closure of the set, which is swept out by all the tangent
lines to $Q$ at the points of the set $A$. Using the estimate
(\ref{0017}) we see that
$\dim J(C,Q)\leq k-1$ and for this reason $J(C,Q)$ is a
proper closed subset of the space ${\mathbb P}$. Furthermore, it is
clear that
$$
T(C,Q)\subset J(C,Q).
$$
Consequently, for a point of general position
$$
a\in{\mathbb P}\setminus J(C,Q)
$$
any line $L_{a,q}$, where $q\in C$, has no other points of
intersection with $Q$ and does not touch $Q$ at the point $q$.
Since the conditions (i) and (ii) are obviously open in $a$ and
$q$, there is an open set $U\subset Q$ that contains the entire
curve $C$ and satisfies (for the fixed $a$) the conditions (i) and
(ii). Q.E.D. for the corollary.

{\bf Corollary 2.} {\it Assume that $\mathop{\rm codim}Q=2$,
$C\subset Q$ is a curve (possibly reducible)
and $q\in C$. For a sufficiently general
choice of the point $a\in{\mathbb P}\setminus Q$ the projection
$$
\pi_a\colon Q\to R\subset{\mathbb P}^{k-1}
$$
is bijective and preserves multiplicity in a neighborhood of the
point $q$. More precisely, there exists an open set $U\subset Q$
such that:

{\rm (i)} $C\cap (Q\setminus U)=\Gamma$ is a finite set of points;

{\rm (ii)} $\pi^{-1}_a(\pi_a(U))=U$, the map $\pi_a\colon U\to
\pi_a(U)$ is bijective and preserves multiplicity;

{\rm (iii)} for each point $z\in\Gamma$ the line $L_{a,z}$ does
not touch $Q$ at $z$, $L_{a,z}\cap Q=\{z,z^*\}$, where $z^*\neq
z$ is a smooth point of the variety $Q$ and $L_{a,z}$ does not
touch $Q$ at the point $z^*$. Moreover, for each $z\in\Gamma$
there exists an open set
$U_z\subset Q$, $z\in U_z$, such that for any point $v\in U_z$
either

\begin{itemize}

\item
the line $L_{a,v}$ meets $Q$ at the point $v$ only and the intersection
is transversal, or

\item
$v$ has the same properties as $z$, that is, the line $L_{a,v}$
meets $Q$ at precisely two points, $v$ and $v^*$, $v^*$ smooth
on $Q$, and in both cases the intersection is transversal.

\end{itemize}
}

{\bf Proof.} In this case $J(C,Q)={\mathbb P}$, so that for any point
$a\in{\mathbb P}$ there are points $q\in C$ where the projection
$\pi_a$ increases the mupltiplicity. However, if $A\subset C$ is
any finite set of points then $\dim J(A,Q)\leq k-1$. Therefore
for a general point $a\in{\mathbb P}$ there exists an open set
$U\subset Q$ such that both conditions (i) and (ii) are satisfied.
(In order to get the condition (i), it is enough for the set $A$
to contain at least one point $a_i\in A$, $a_i\in C_i$, for each
irreducible component $C_i$ of the curve $C$.)

Furthermore, $\dim (C,\mathop{\rm Sing}Q)\leq k-1$, so that for
each point $z\in\Gamma$ the line $L_{a,z}$ meets $Q$ (leaving
aside the very point $z$) at smooth points only. Counting
dimensions, it is easy to see that for a general point $a$ the
intersection $L_{a,z}$ consists of precisely two points $z$,
$z^*$, whereas at $z^*$ this intersection is transversal.
Finally, it is easy to see that the property of the line
$L_{a,v}$ meeting $Q$ transversally at at most two points
($v$ itself and possibly another one, $v^*$, which is smooth
on $Q$) is open in $v$. This proves our last claim. Q.E.D.

Now it is trivial to generalize Corollaries 1 and 2 for effective
algebraic cycles (instead of irreducible subvarieties). Let
$$
Q=\sum_{i\in I_Q}m_iQ_i
$$
be an effective cycle with $Q_i$ pair-wise distinct and of the
same dimension $\mathop{\rm dim}Q_i=\mathop{\rm dim}Q$, $m_i\geq 1$.
We define the multiplicity of the cycle $Q$ at a point
$y\in{\mathbb P}$ by linearity:
$$
\mathop{\rm mult}\nolimits_yQ=
\sum_{i\in I_Q}m_i\mathop{\rm mult}\nolimits_yQ_i.
$$
It is easy to see that the claim of Corollary 1 remains true if
we replace an irreducible subvariety $Q$ by an effective cycle $Q$,
and $R=\pi_a(Q)$ by
$$
R=(\pi_a)_*Q=\sum_{i\in I_Q}m_iR_i,
$$
where $R_i=\pi_a(Q_i)$ and we may assume that the subvarieties
$R_i$ are again pair-wise distinct.

The claims (i) and (ii) of Corollary 2 also remain true after the
same modification. As for the claim (iii), it should read as
follows:

(iii) {\it for each point $z\in \Gamma$ the line $L_{a,z}$ does not
touch $\mathop{\rm Supp}Q$ at $z$,
$L_{a,z}\cap \mathop{\rm Supp}Q=\{z,z^*\}$, where $z^*\neq z$ is a
smooth point of the closed set $\mathop{\rm Supp}Q$, and $L_{a,z}$
does not touch $\mathop{\rm Supp}Q$ at the point $z^*$. In other
words, there exists one and only one component $Q_{i(z)}$ of the
cycle $Q$, $i(z)\in I_Q$, such that $Q_{i(z)}\ni z^*$ and
$Q_{i(z)}$ is smooth at $z^*$. In particular,

\begin{equation}
\label{0019}
\mathop{\rm mult}\nolimits_{\pi_a(z)}R=
\mathop{\rm mult}\nolimits_zQ+m_{i(z)}\leq
\mathop{\rm mult}\nolimits_zQ+\mathop{\rm max}\limits_{i\in I_Q}\{m_i\}.
\end{equation}

Moreover, for each point $z\in \Gamma$ there exists an open set
$U_z\subset\mathop{\rm Supp}Q$,
$z\in U_z$, such that for any point $v\in U_z$ either

\begin{itemize}

\item
the line $L_{a,v}$ meets $\mathop{\rm Supp}Q$ at the point $v$ only
and the intersection is transversal, or

\item
$v$ has the same properties as $z$, that is, the line $L_{a,v}$
meets $\mathop{\rm Supp}Q$ at precisely two points, $v$ and $v^*$,
$v^*$ smooth on $\mathop{\rm Supp}Q$,
and in both cases the intersection is transversal.

\end{itemize}

In particular, for any $v\in U_z$ we have the estimate
\begin{equation}
\label{0023}
\mathop{\rm mult}\nolimits_{\pi_a(v)}R
\leq \mathop{\rm mult}\nolimits_vQ+\mathop{\rm max}\limits_{i\in I_Q}\{m_i\}.
\end{equation}
}

{\bf Proof of the claim (iii).} All the assertions of the claim but
the estimate (\ref{0019}) follow directly from part (iii) of
Corollary 2. So let us prove this estimate. For a general 2-plane
$P\subset{\mathbb P}^{k-1}$, containing the point $r=\pi_a(z)$, we
have
$$
\mathop{\rm mult}\nolimits_rR=(P\cdot R)_r.
$$
As in the proof of Proposition 6 above, we take the inverse image
$P^+=\overline{\pi_a^{-1}(P)}$ of the plane $P$ and obtain the
equality
\begin{equation}
\label{0021}
(P\cdot R)_r=(P^+\cdot Q)_z+(P^+\cdot Q)_{z^*}.
\end{equation}
Since the line $L_{a,z}$ does not touch $\mathop{\rm Supp}Q$ at $z$
and $z^*$, and $z^*$ is a smooth point of $Q_{i(z)}$,
we may rewrite (\ref{0021}) as the left-hand side of the
inequality (\ref{0019}). The rest is obvious. Q.E.D.

In what follows we shall refer to Corollaries 1 and 2 in their
generalized form just discussed (that is, for effective equidimensional
cycles $Q$).

\subsection{Singularities of the divisor $F$}

Let us prove part (i) of Proposition 4. Note that if a subvariety
$W\subset X$ of a smooth variety $X$ is a log-centre of the
pair $(X,D)$, where $D=\sum\limits_{i\in I_D}d_iD_i$ is a
${\mathbb Q}$-divisor, then for each point $x\in W$
$$
\mathop{\rm mult}\nolimits_x D=\sum\limits_{i\in I_D}
d_i \mathop{\rm mult}\nolimits_x D_i>1.
$$
Therefore if a point $x\in{\mathbb P}^{M-2}$ lies on a log-centre
of the pair $({\Bbb P}^{M-2},\frac{1}{2n^2}F)$ then
$$
\mathop{\rm mult}\nolimits_x F >2n^2.
$$

Denote by $C\ni x$ the connected component of the set
$$
\{z\in Z\,|\, \mathop{\rm mult}\nolimits_x Z>n^2\},
$$
containing the point $x$. Since $Z\sim n^2H^2$ is an effective cycle
of codimension two on the smooth hypersurface $V$, we conclude
that by Proposition 5 that $C$ is either a point or a curve.

Let us represent the projection from a line $L$ as a composition
of two projections: $\pi_L=\pi_a\circ \pi_b$, where $b\in{\Bbb
P}$ and $a\in{\mathbb P}^{M-1}$ are general points. Set
$Z_b=(\pi_b)_*Z$. By Corollary 1 in
$\mathop{\rm Supp}Z$ there exists an open subset $U_b\supset C$ such that
$\pi_b|_{U_b}$ preserves multiplicity. Consequently, the connected
closed algebraic set $C_b=\pi_b(C)$ is a connected component of the set
$$
\{z\in \mathop{\rm Supp}Z_b\, |\, \mathop{\rm mult}\nolimits_z Z_b>n^2\}.
$$
If $C=\{x\}$ is a point, then $C_b=\pi_b(x)$ is again a point and we
apply Corollary 1 once again to conclude that for any point $z$
in a neighborhood of the point $y$ the preimage
$(\pi_L|_{\mathop{\rm Supp}Z})^{-1}(z)$ consists of only one point and if
$z\neq y=\pi_a(\pi_b(x))=\pi_L(x)$ then we have
$$
\mathop{\rm mult}\nolimits_z F=\mathop{\rm
mult}\nolimits_{(\pi_L|_{\mathop{\rm Supp}Z})^{-1}(z)}Z\leq n^2,
$$
so that the pair $({\mathbb P}^{M-2},\frac{1}{2n^2}F)$ is
log-terminal at the point $y$ in dimension one.

So we assume that $C$ is a connected curve. Therefore
$C_b\subset \mathop{\rm Supp}Z_b$ is also a connected curve and we
apply Corollary 2 (in the generalized form) to the effective
cycle $Z_b$ in a neighborhood of the curve $C_b$. We see
that there
exists an open subset $U_a\subset\mathop{\rm Supp}Z_b$, containing the point
$\pi_b(x)$ and meeting each component of the curve $C_b$, such
that the map
$\pi_a\colon U_a\to \pi_a(U_a)\subset\pi_L(\mathop{\rm Supp}Z)$ is
bijective and preserves multiplicity. In particular,
$$
\pi_L(C)\cap\pi_a(U_a)=\{z\in\pi_a(U_a)\,|\, \mathop{\rm
mult}\nolimits_z F>n^2\}.
$$
Therefore any variety $W\subset \mathop{\rm Supp}F$ which has a
non-empty intersection with $\pi_a(U_a)$ and is not contained
entirely in $\pi_L(C)$ cannot be a log-centre of the pair $({\Bbb
P}^{M-2},\frac{1}{2n^2}F)$.

It remains to consider any of the ``bad'' points $p\in\Gamma=
C_b\cap(\mathop{\rm Supp}Z_b\setminus U_a)$. The line $L_{a,p}$ is not tangent to
$Z_b$ at the point $p$ and
$$
L_{a,p}\cap \mathop{\rm Supp}Z_b=\{p,p^*\},
$$
where $p^*\neq p$ is a smooth point of one of the irreducible
components of $Z_b$, and moreover the intersection of $L_{a,p}$ and
$\mathop{\rm Supp}Z_b$ at the point $p^*$ is transversal.
By Corollary 2, part (iii) in the generalized form we get a
neighborhood $U_p\ni p$ such that for any $v\in U_p$ the estimate
(\ref{0023}) holds. Since the multiplicity of any irreducible component
of the cycle $Z$ is not higher than $n^2$ (by the Lefschetz theorem)
and the same is accordingly true for $Z_b$, we get finally for any
point $q\in U_p$:
\begin{equation}
\label{3}
\mathop{\rm mult}\nolimits_{\pi_a(q)}F
\leq \mathop{\rm mult}\nolimits_q Z_b+n^2.
\end{equation}

Now recall that $C_b$ is a connected component of the set $\{z\in
Z_b\,|\, \mathop{\rm mult}\nolimits_z Z_b>n^2\}$. By (\ref{3})
this implies that in a neighborhood of the point $\pi_a(p)$ the
connected set $\pi_a(C_b)$ contains all the points $z\in F$ where
$$
\mathop{\rm mult}\nolimits_z F\geq 2n^2+1.
$$
Consequently, the curve $\pi_a(C_b)$ has a neighborhood $U\subset
\mathop{\rm Supp}F$, $U\supset \pi_a(C_b)$ such that
$$
\pi_a(C_b)\supset \{z\in U\, |\, \mathop{\rm mult}\nolimits_z
F\geq 2n^2+1\}.
$$
Therefore the pair $({\mathbb P}^{M-2},\frac{1}{2n^2}F)$ is
log-canonical at the point $y=\pi_L(x)$ in dimension 2.

\section{Log-canonical singularities}

In this section we prove Proposition 1. The arguments are based
on the connectedness principle of Shokurov and Koll\' ar.

\subsection{The connectedness principle of Shokurov and Kollar}

We will need a very particular case of this principle. Let $X$ be
a smooth variety, $D=\sum\limits_{i\in I_D}d_i D_i$ a ${\Bbb
Q}$-divisor, $d_i>0$, where the prime divisors $D_i$ are
pair-wise distinct. Let
$$
f\colon Y\to X
$$
be a resolution of singularities, $\{E_i\,|\, i\in I_f\}$ the set
of exceptional divisors, where the divisor
$$
\mathop{\bigcup}\limits_{i\in I_D}\widetilde D_i\,\cup\,
\mathop{\bigcup}\limits_{i\in I_f}E_i
$$
has normal crossings on $Y$, $\widetilde D_i$ stands for the
strict transform of $D_i$ on $Y$. Set
$$
\widetilde D=\sum_{i\in I_D}d_i\widetilde D_i
$$
and write down
\begin{equation}
\label{0029}
K_Y+\widetilde D=f^*(K_X+D)+\sum_{i\in I_f}e_iE_i.
\end{equation}

{\bf The connectedness principle (Shokurov [Sh], Kollar [K]).}
{\it Assume that the class $-(K_X+D)$ is numerically effective
and big. Then the closed algebraic set
$$
\mathop{\bigcup}\limits_{d_i\geq 1} D_i \cup
f(\mathop{\bigcup}\limits_{e_i\leq -1}E_i)
$$
is connected.}

{\bf Proof} (extracted from [K], 17.4) is based on the

{\bf Kawamata-Viehweg vanishing theorem.} {\it Let $Y$ be a smooth
projective variety, $R=\sum\limits_{i\in I_R}r_iR_i$ an effective
${\mathbb Q}$-divisor with $0<r_i<1$, $R_i$ pair-wise distinct
and
$$
\mathop{\bigcup}\limits_{i\in I_R} R_i
$$
with normal crossings. Assume that $L\in\mathop{\rm Pic}Y$ is a
class such that the class $(L-R)$ is numerically effective and big,
that is, $(L-R)^{\dim Y}>0$. Then
$$
H^j(Y,K_Y+L)=0
$$
for} $j\geq 1$.

For a proof, generalizations and explanations see [Kw,V,EV].

Now, following Koll\' ar, let us obtain the connectedness
principle. For convenience of notations set
$$
J=I_f\cup I_D
$$
and for $i\in I_D$ set $E_i={\widetilde D}_i$, $e_i=-d_i$.
Now we can rewrite (\ref{0029}) as
\begin{equation}
\label{0031}
K_Y=f^*(K_X+D)+\sum_{i\in J}e_iE_i.
\end{equation}
For each $i\in J$ set
$$
m_i=\min\{m\in{\mathbb Z}|m\geq e_i\},\quad \alpha_i=m_i-e_i\geq 0,
$$
so that $e_i=m_i-\alpha_i$, $\alpha_i<1$. Rewrite (\ref{0031}) as
$$
-f^*(K_X+D)=-K_Y+\sum_{i\in J}m_iE_i-\sum_{i\in J}\alpha_iE_i.
$$
Since by assumption $-(K_X+D)$ is nef and big, the same is true
for its $f$-pull back. Now apply the vanishing theorem to
$$
L=-K_Y+\sum_{i\in J}m_iE_i \quad\mbox{and}\quad
R=\sum_{i\in J}\alpha_iE_i.
$$
We obtain the vanishing
$$
H^j(Y,\sum_{i\in J}m_iE_i)=0
$$
for $j\geq 1$. On the other hand $m_i\in{\mathbb Z}_+$ if and only
if $e_i>-1$ and moreover $m_i\geq 1$ if and only if $e_i>0$. Therefore
setting
$$
E^+=\sum_{e_i>-1}m_iE_i=\sum_{e_i>0}m_iE_i
$$
and
$$
E^-=-\sum_{e_i\leq -1}m_iE_i,
$$
we see that
\begin{equation}
\label{0033}
H^j(Y,E^+-E^-)=0
\end{equation}
for $j\geq 1$ where both $E^+$, $E^-$ are effective divisors. Now apply
(\ref{0033}) to the standard exact sequence
$$
0\longrightarrow {\cal O}_Y(E^+-E^-)
\longrightarrow {\cal O}_Y(E^+)
\longrightarrow {\cal O}_Y(E^+)|_{E^-}
\longrightarrow 0.
$$
We conclude that the map
\begin{equation}
\label{0035}
H^0(Y,E^+)\longrightarrow H^0(E^-,E^+|_{E^-})
\end{equation}
is surjective. But since for $i\in I_D$ all the coefficients $e_i$
are negative, the divisor $E^+$ is $f$-exceptional. Thus
$$
H^0(Y,E^+)\cong{\mathbb C}.
$$
On the other hand, the components $E_i$ of the divisors $E^+$ and $E^-$
form two disjoint sets, so that the restriction $E^+|_{E^-}$ is an
effective divisor on the scheme $E^-$.

Let $l\geq 1$ be the number of connected components of the set
$$
\mathop{\bigcup}\limits_{e_i\leq -1}E_i
$$
which is precisely $\mathop{\rm Supp}E^-$. Obviously,
$$
\dim H^0(E^-,E^+|_{E^-})\geq l.
$$
Therefore (\ref{0035}) gives a surjective map
$$
{\mathbb C}\longrightarrow{\mathbb C}^l\longrightarrow 0
$$
which of course implies that $l\leq 1$. In other words,
$\mathop{\rm Supp}E^-$ is connected, which is precisely
what we need. Q.E.D. for the connectedness principle.

\subsection{ Log-canonical singularities of divisors in ${\mathbb P}^k$}

Let us prove Proposition 1. It is sufficient to prove the
``terminal'' version of the proposition, since the ``canonical''
version is its direct implication. Moreover, since the property
of being log-terminal is open with respect to the coefficients
$d_i$, we may assume that the inequality (\ref{1a}) is strict,
that is,
$$
l+\mathop{\rm deg} D < k+1.
$$
Now assume the converse:
$$
x\in\mathop{\rm LC} (X,D).
$$
By assumption $\dim\mathop{\rm LC} (X,D,x)\leq l-1$, so that for
a general plane $P\subset{\mathbb P}^k$ of codimension $l$ we get
$$
P\cap \mathop{\rm LC} (X,D,x)=\emptyset.
$$
To simplify our notations, we write in this section ${\mathbb P}$
instead of ${\mathbb P}^k$.

Fix such a plane $P$. Denote by $\Lambda_P\subset |H|$ the linear
system of hyperplanes containing $P$. If $l=1$, then
$\Lambda_P=\{P\}$. In this case set $H^{\sharp}=P$.

If $l\geq 2$, then the linear system $\Lambda_P$ is free from
fixed components, $\mathop{\rm Bs} \Lambda_P=P$. In this case set
$$
H^{\sharp}=\sum_{i\in I_H}\varepsilon H_i,
$$
where $\varepsilon>0$ is a small rational number, $\{H_i\, |\,
i\in I_H\}$ is a general set of hyperplanes in $\Lambda_P$ and
$$
\mathop{\rm deg} H^{\sharp}=l=\varepsilon\cdot\sharp I_H.
$$
Consider a new ${\Bbb Q}$-divisor $D^{\sharp}=D+H^{\sharp}$.

{\bf Lemma 5.} {\it Outside the plane $P$ the closed sets
$$
\mathop{\rm LC} ({\mathbb P},D)\quad\mbox{and}\quad \mathop{\rm
LC}({\mathbb P},D^{\sharp})
$$
coincide.}

{\bf Proof.} This is almost obvious. Let $f\colon Y\to{\mathbb P}$
be a resolution of the pair $({\mathbb P},D)$. Clearly outside $P$
the divisor $H^{\sharp}$ has normal crossings and
$$
\widetilde H^{\sharp}|_{Y\setminus f^{-1}(P)}=
f^{-1}(H^{\sharp})|_{Y\setminus f^{-1}(P)},
$$
so that on $Y\setminus f^{-1}(P)$ we get
$$
K_Y+\widetilde D^{\sharp}= f^*(K_{{\mathbb P}}+D^{\sharp})+\sum_{i\in
I_f}e_iE_i.
$$
Moreover, since on $Y\setminus f^{-1}(P)$ the linear system
$f^*\Lambda_P$ is free, the divisor
$$
\mathop{\bigcup}\limits_{i\in I_D}\widetilde D_i\,\cup\,
\mathop{\bigcup}\limits_{i\in I_f} E_i\,\cup\,
\mathop{\bigcup}\limits_{i\in I_H}\widetilde H_i
$$
has on $Y\setminus f^{-1}(P)$ normal crossings. Taking into
account that $\varepsilon>0$ is small, we get finally:
$$
\mathop{\rm LC} ({\mathbb P}\setminus P, D^{\sharp})=
f\left(\mathop{\bigcup}\limits_{e_i\leq -1}
E_i\right)=\mathop{\rm LC} ({\mathbb P}\setminus P, D),
$$
which is what we need.

Now write down
$$
\begin{array}{c}
\displaystyle K_Y=f^*(K_{{\mathbb P}}+D^{\sharp})-\sum_{i\in
I_D}d_i\widetilde D_i-\sum_{i\in I_H}\varepsilon\widetilde H_i+
\\  \\
\displaystyle +\sum_{i\in I^*_f}e_iE_i+e_PE
\end{array}
$$
for some resolution $f\colon Y\to {\mathbb P}$ of the pair $({\Bbb
P},D^{\sharp})$, where for $l=1$ we have $\varepsilon=0$,
$E=\widetilde P$ and $e_P=-1$, and for $l\geq 2$ the exceptional
divisor $E$ is determined by the condition that $f(E)=P$ and in a
neighborhood of the general point of $E$ the morphism $f\colon
Y\to{\mathbb P}$ is the blow up of the plane $P$.

{\bf Lemma 6.} $e_P=-1$.

{\bf Proof:} indeed,
$$
\sum_{i\in I_H}\varepsilon =l,
$$
whereas the discrepancy of $E$ is equal to $a(E)=\mathop{\rm
codim} P-1=l-1$, whence we get the claim of the lemma.

Note that $-(K_{{\mathbb P}}+D^{\sharp})$ is a numerically effective
and big class: more precisely, it is equal to $(k+1-l-\mathop{\rm
deg} D)H$, where the coefficient in brackets is strictly positive.
Now by the connectedness principle and the two previous lemmas
the set
$$
f\left(\mathop{\bigcup}\limits_{e_i\leq -1}E_i\right)\cup
P=\mathop{\rm LC} ({\mathbb P},D)\cup P
$$
is connected. However, the connected component $\mathop{\rm
LC}({\mathbb P},D,x)$ by assumption is non-empty and moreover
$$
\mathop{\rm LC}({\mathbb P},D,x)\cap P=\emptyset.
$$
Therefore the closed set $\mathop{\rm LC}({\mathbb P},D)\cup P$ has
at least two connected components. A contradiction. Q.E.D. for
Proposition 1.

\section{The main construction}
In this section we prove part (ii) of Proposition 4. The argument
is based on the geometric construction introduced in [P2].

\subsection{The direct image of a maximal singularity}
Let
$$
\pi\colon{\mathbb P}-\,-\,\to{\mathbb P}^{M-2}
$$
be the linear projection from a general line $L$,
$$
\pi_L|_V\colon
V-\,-\,\to{\mathbb P}^{M-2}
$$
its restriction onto the hypersurface
$V$. The map $\pi_L|_V$ is not well defined at the $M$ points of
intersection $L\cap V$. Blow these $M$ points up:
$$
\sigma\colon\widetilde V\to V
$$
and denote the regular extension of the map $\pi_L|_V$ on
$\widetilde V$ by the symbol
$$
\pi\colon\widetilde V\to{\mathbb P}^{M-2}.
$$
The following properties are satisfied when the line $L$ is
sufficiently general:

(1) $L\cap \mathop{\rm Supp}Z=\emptyset$,
so that for the strict transform
$\widetilde Z$ of the effective cycle $Z=(D_1\circ D_2)$ on
$\widetilde V$ we have $\widetilde Z=\sigma^{-1}(Z)$;

(2) the curve $\pi^{-1}(y)$ is smooth
(recall that $y=\pi_L(x)$), meets $Z$ at the unique
point $x$ and is not tangent to $Z$ at $x$;

(3) the morphism $\pi$ is birational on each irreducible
component of the cycle $Z$.

Let us consider the direct image $\pi_*\nu$ of the discrete
valuation $\nu$. Algebraically $\pi_*\nu$ can be defined (up to
dividing by an integer) as the composition
$$
\pi_*\nu \colon {\mathbb C}({\mathbb P}^{M-2})^{\times}
\stackrel{\pi^*}{\hookrightarrow} {\mathbb C}(\widetilde V)^{\times}=
{\mathbb C}(V)^{\times}\stackrel{\nu}{\longrightarrow}
{\mathbb Z},
$$
where ${}^{\times}$ means taking the multiplicative group of
non-zero elements of a field. Realizing $\pi_*\nu$ geometrically, let
$$
\varphi\colon X\to{\mathbb P}^{M-2}
$$
be a birational morphism with $X$ smooth, the set
$\{E_i\,|\,i\in I_X\}$ of exceptional divisors having normal crossings
and
\begin{equation}
\label{4}
\varphi\left(\mathop{\bigcup}\limits_{i\in I_X}
E_i\right)=y,
\end{equation}
such that for some exceptional divisor $E=E_a$ we get
$$
\pi_*\nu=\mathop{\rm ord}\nolimits_E(\cdot).
$$
Geometrically, it means the following. Take
$$
V^*=\widetilde V\mathop{\times}\nolimits_{{\mathbb P}^{M-2}}X.
$$
Taking into account (\ref{4}) and the fact that the fiber
$\pi^{-1}(y)$ is non-singular, we see that the variety $V^*$ is
smooth and the projection
$$
f\colon V^*\to\widetilde V
$$
is a birational morphism with the set of exceptional divisors
$\{E^*_i\,|\,i\in I_X\}$, and moreover
$$
E^*_i=C\times E_i\quad\mbox{and}\quad f
\left(\mathop{\bigcup}\limits_{i\in I_X} E^*_i\right)=C,
$$
where $f|_{E^*_i}$ is the projection onto the first factor.
Moreover, the following diagram is commutative
$$
\begin{array}{rcccl}
   &   V^*  &  \stackrel{f}{\longrightarrow} & \widetilde V &   \\
\pi &  \downarrow  &     &    \downarrow &  \pi \\
   &   X  &   \stackrel{\varphi}{\longrightarrow} & {\mathbb P}^{M-2},
\end{array}
$$
where, somewhat abusing our notations, we denote the natural
projection of $V^*$ onto $X$ also by the symbol $\pi$, since
outside $\cup E^*_i$ the variety $V^*$ is isomorphic to
$\widetilde V$ and the left-hand projection extends the
right-hand one.

Now by the definition of direct image we get
$$
\pi(\mathop{\rm centre}(\nu,V^*))=E.
$$
Since $\mathop{\rm centre}(\nu,\widetilde V)=x$, we can say more.
For each $i\in I_X$ set
$$
\Delta_i=E^*_i\cap f^{-1}(x)\subset E^*_i.
$$
Obviously all the projections $\pi\colon\Delta_i\to E_i$ are
isomorphisms. Now we get
$$
\mathop{\rm centre}(\nu,V^*)=\Delta,
$$
where $\Delta=\Delta_a\subset E^*=E^*_a$.

Let $\Sigma^*$ be the strict transform of the linear system
$\Sigma$ on $V^*$.

{\bf Lemma 7.} {\it The strict transform coincides with the pull
back:}
$$
\Sigma^*=f^*\widetilde \Sigma=f^*\sigma^*\Sigma.
$$

{\bf Proof.} As it was mentioned above, for each exceptional
divisor $E^*_i$ we have $f(E^*_i)=C$, but the curve $C$ is not a
base curve for $\widetilde \Sigma$, since
$$
C\not\subset\widetilde Z.
$$
Q.E.D. for the lemma.

\subsection{Comparing multiplicities}
Here we use the ring structure on the group of classes of cycles
on smooth varieties, the operations of direct image and flat
inverse image described in [F].

Take the divisors
$$
D^*_i=f^*\widetilde D_i=f^*\sigma^* D_i,
$$
$i=1,2$, and let
$$
Z^*=(D^*_1\circ D^*_2), \quad
\widetilde Z=(\widetilde D_1\circ \widetilde D_2)
$$
be the algebraic cycles of scheme-theoretic
intersection of these divisors. Denote by $z^*$, $\tilde z$ and
$z$ the classes of $Z^*$, $\widetilde Z$ and $Z$ in $A^2 V^*$,
$A^2\widetilde V$ and $A^2 V$, respectively.

{\bf Corollary 3 (from Lemma 7).} {\it The following equality holds:
$z^*=f^*\tilde z=f^*\sigma^* z$.}

{\bf Proof.} Indeed, $f^*$ and $f^*\sigma^*$ are ring homomorphisms.
Q.E.D.

Now let us compare the direct image $F$ of the cycle $Z$ on
${\mathbb P}^{M-2}$ and the direct image of $Z^*$ on $X$.

{\bf Lemma 8.} {\it The following equality is true: $\pi_*
Z^*=\varphi^* F$.}

{\bf Proof.} By construction $\pi_*Z=F$; therefore,
$\pi_*\widetilde Z=F$. Thus we get:
$$
\pi_* Z^*=\widetilde F+\sum_{i\in I_X}b_iE_i,
$$
where $\widetilde F$ is the strict transform of the cycle $F$ on
$X$, $b_i\in{\mathbb Z}$ are uniquely determined integers. On the
other hand,
$$
\varphi^* F=\widetilde F+\sum_{i\in I_X}c_iE_i,
$$
where $c_i=\mathop{\rm ord}\nolimits_{E_i}F\in{\mathbb Z}_+$.
However, by the previous corollary (and the standard theorems of
intersection theory [F, 1.7 and Ch. 8])
$$
\pi_* z^*=\pi_*(f^*\sigma^* z)=\varphi^*\pi_* \tilde z
$$
(here $f^*$, $\varphi^*$ and $\sigma^*$ are pull back operations
in $A^*$, $\pi_*$ is the direct image, preserving dimension), so
that the classes of divisors
$$
\sum_{i\in I_X}b_iE_i\quad\mbox{and}\quad \sum_{i\in I_X}c_iE_i
$$
in $A^1X$ coincide. Consequently, $b_i=c_i$ for all $i\in I_X$,
which proves the lemma.

{\bf Corollary 4.} {\it The following equality holds:}
\begin{equation}
\label{5}
\nu_E(F)=\mathop{\rm mult}\nolimits_{\Delta}Z^*.
\end{equation}

Here the right-hand side means the multiplicity with which the
irreducible subvariety $\Delta\subset V^*$ of codimension two
comes into the cycle $Z^*$.

\subsection{Computation of $\nu_E(F)$}
The right-hand side of the equality (\ref{5}) is not hard to
compute (actually, this has already been done in [P2]). Namely,
let
$$
\begin{array}{rccc}
\gamma_{i,i-1}\colon & V^{\sharp}_i & \longrightarrow &
V^{\sharp}_{i-1} \\
    &  \bigcup &   &   \bigcup \\
    &  E^{\sharp}_{i} & \longrightarrow & B_{i-1}
\end{array}
$$
be the sequence of blow ups of the centres of the valuation $\nu$,
starting with $V^{\sharp}_0=V^*$, $B_0=\Delta$, that is,
$$
B_i=\mathop{\rm centre}(\nu,V^{\sharp}_i),
$$
$E^{\sharp}_i=\gamma^{-1}_{i,i-1}(B_{i-1})$, $i=1,\dots,N$,
$E^{\sharp}_N\subset V^{\sharp}_N$ is a divisorial realization of the
valuation $\nu$. Let $\Sigma^i$ be the strict transform of the
linear system $\Sigma$ on $V^{\sharp}_i$. Define the
multiplicities $\mu_i$ as
$$
\mu_i=\mathop{\rm mult}\nolimits_{B_{i-1}}\Sigma^{i-1}.
$$
In the usual way (see [IM,IP,P1-P8]) we introduce the graph
$\Gamma$ of exceptional divisors $E^{\sharp}_i$, $i=1,\dots,N$,
taking the set
$$
\{E^{\sharp}_i\}\leftrightarrow \{i=1,\dots,N\}
$$
to be the set of vertices, and drawing an arrow
$$
i\longrightarrow j
$$
if and only if $i>j$ and the subvariety $B_{i-1}$ lies in the
strict transform $(E^{\sharp}_j)^{i-1}$ of the exceptional
divisor $E^{\sharp}_j$ on $V^{\sharp}_{i-1}$.
Let the integers $p_{ij}$ for $i>j$ stand for the number of paths
$$
i_0=i\longrightarrow i_1 \longrightarrow \dots
\longrightarrow i_k=j,
$$
$k=1,\dots$, from $i$ to $j$ in the graph $\Gamma$ of the
valuation $\nu$ and set $p_{ii}=1$, so that
$$
p_{ij}=\nu_{E^{\sharp}_i}(E^{\sharp}_j).
$$
By the standard formulae of the intersection theory [P3, Appendix
B], see also [P4], we get
$$
\mathop{\rm mult}\nolimits_{\Delta}Z^*\geq \sum^N_{i=1}\mu_i^2.
$$

{\bf Lemma 9.} {\it The following estimate is true:
\begin{equation}
\label{5a}
\sum^N_{i=1}p_i\mu_i>n\cdot\left(
a(E)p_1+\sum^N_{i=1}p_i\right)
\end{equation}
where $p_i=p_{N,i}$ and $a(E)$ stands for the discrepancy of} $E$.

{\bf Proof.} By the Noether-Fano inequality
$$
\nu(\Sigma)>n\cdot \mathop{\rm discrepancy}(\nu).
$$
Here $\nu(\Sigma)=\nu(\Sigma^*)$, since $\Sigma^*$ is the pull
back of the linear system $\Sigma$ on $V^*$. By the standard
technique (see [P2,P3,P5])
$\nu(\Sigma^*)$ is the left-hand side of (\ref{5a}).
Let us estimate the discrepancy of the valuation $\nu$ {\it with
respect to} $V$: we get
$$
K_{V^*}=f^*\sigma^* K_V+\sum_{i\in I_X}a(E_i)E^*_i,
$$
so that
$$
a(\nu,V)\geq a(E,{\mathbb P}^{M-2})\cdot \nu(E^*)+a(\nu,V^*).
$$
Since $\mathop{\rm centre}(\nu,V^*)=\Delta\subset E^*$, we obtain
the inequality
$$
\nu(E^*)\geq \nu(E^{\sharp}_1)=p_1
$$
(in fact, one can say more precisely that
$$
\nu(E^*)=\sum_{B_{i-1}\subset (E^*)^{i-1}}p_i,
$$
where $(E^*)^{i-1}$ is the strict transform of $E^*$ on
$V^{\sharp}_{i-1}$, but we do not need that). Finally, all the
centres $B_{i-1}$ are of codimension 2, so that
$$
a(\nu,V^*)=\sum^N_{i=1}p_i,
$$
which is what we need. Q.E.D. for Lemma 9.

{\bf Corollary 5.} {\it The following estimate holds:}
\begin{equation}
\label{6}
\nu_E(F)>n^2\cdot\frac{\displaystyle\left(a(E)p_1+\sum^N_{i=1}p_i\right)^2}{
\displaystyle\sum^N_{i=1}p^2_i}.
\end{equation}

{\bf Proof.} Compute the minimum of the quadratic form
$$
\sum^N_{i=1}\mu^2_i
$$
on the hyperplane
\begin{equation}
\label{8} \sum^N_{i=1}p_i\mu_i=C.
\end{equation}
The minimum is attained for $\mu_i=p_i\lambda$, where the common
constant $\lambda$ can be obtained from (\ref{8}). Elementary
computations complete the proof.

{\bf Corollary 6.} {\it The following estimate holds:}
\begin{equation}
\label{9} \nu_E(F)>2n^2(a(E)+1).
\end{equation}

{\bf Proof.} Assume first that $N\geq 2$. Opening the brackets in
(\ref{6}), we obtain the inequality
$$
\left(a(E)p_1+\sum^N_{i=1}p_i\right)^2>
2a(E)p_1\left(\sum^N_{i=1}p_i\right)+\left(\sum^N_{i=1}p_i\right)^2.
$$
But $p_1\geq p_i$ for all $i=1,\dots,N$, so that $p_1p_i\geq
p_i^2$ and the first component of the right-hand side is not
smaller than
$$
2a(E)\sum^N_{i=1}p^2_i.
$$
By the definition of the integers $p_i$ we get
$$
p_1=\sum_{i\to 1}p_i,
$$
so that
$$
\begin{array}{c}
\displaystyle \left(\sum^N_{i=1}p_i\right)^2\geq
p_1\cdot\sum^N_{i=1}p_i+\left(\sum_{i\to
1}p_i\right)\sum^N_{i=1}p_i= \\   \\
\displaystyle =2p_1\sum^N_{i=1}p_i\geq 2\sum^N_{i=1}p^2_i,
\end{array}
$$
whence by (\ref{6}) we get the estimate (\ref{9}) for $N\geq 2$.

If $N=1$, then $p_1=1$ and
$$
\nu_E(F)>n^2(a(E)+1)^2,
$$
but $(a(E)+1)^2=2a(E)+a(E)^2+1\geq 2a(E)+2$, which is what we
need.

It remains to note that the estimate (\ref{9}) means precisely
that the pair $({\mathbb P}^{M-2},\frac{1}{2n^2}F)$ is not
log-canonical at the point $y$.

The proof of Proposition 4 and thus of our theorem is complete.
\vspace{1cm}

\section*{References}

\noindent [IM] V.A. Iskovskikh and Yu.I. Manin, Three-dimensional
quartics and counterexamples to the L{\" u}roth problem, Math.
USSR Sb., {\bf 15} (1971), no. 1, 141-166. \vspace{0.5cm}

\noindent [C] A. Corti, Singularities of linear systems and
3-fold birational geometry,  In: ``Explicit Birational Geometry of
Threefolds'', London Mathematical Society Lecture Notes Series
{\bf 281}, Cambridge University Press 2000, 259-312. \vspace{0.5cm}

\noindent [CM] A. Corti and M. Mella, Birational geometry of
terminal quartic 3-folds. I, preprint, arXiv:math.AG/0102096
\vspace{0.5cm}

\noindent [CPR] A. Corti, A. Pukhlikov and M. Reid, Fano 3-fold
hypersurfaces, In: ``Explicit Birational Geometry of
Threefolds'', London Mathematical Society Lecture Notes Series
{\bf 281}, Cambridge University Press 2000, 175-258. \vspace{0.5cm}

\noindent [Ch1] I.A. Cheltsov, On a smooth four-dimensional quintic,
Sbornik: Mathematics {\bf 191} (2000), no.9-10, 1399-1419. \vspace{0.5cm}

\noindent [Ch2] I.A. Cheltsov, On the sextic, septic and octic,
preprint, 2000. \vspace{0.5cm}

\noindent [Ch3] I.A. Cheltsov, Log-canonical thresholds on
hypersurfaces, Sbornik: Mathematics {\bf 192} (2001), no. 7-8,
1241-1257. \vspace{0.5cm}

\noindent [ChPk] I.A. Cheltsov and J. Park, Log-canonical
thresholds and generalized Eckard points, preprint,
arXiv:math.AG/0003121. \vspace{0.5cm}

\noindent [EV] H. Esnault and E. Viehweg, Lectures on vanishing
theorems, DMV-Seminar. Bd. 20. Birkh\" auser, 1992. \vspace{0.5cm}

\noindent [F] W. Fulton, Intersection theory. Springer,
1984.\vspace{0.5cm}

\noindent [G1] M. Gizatullin, Fano's inequality
is a mistake, preprint,
arXiv:math.AG/0202069 (2002).
\vspace{0.5cm}

\noindent [G2] M. Gizatullin, Fano's inequality
is also false for three-dimensional
quadric, preprint, arXiv:math.AG/0202117 (2002).
\vspace{0.5cm}

\noindent [G3] M. Gizatullin, Fano's inequality is false
for a simple Cremona
transformation of five-dimensional projective space, preprint,
arXiv:math.AG/0202138 (2002).
\vspace{0.5cm}

\noindent [I] V.A. Iskovskikh, Birational rigidity of Fano
hypersurfaces from the viewpoint of Mori theory, Russian Math.
Surveys {\bf 56} (2001), no. 2, 3-86.\vspace{0.5cm}

\noindent [IP] V.A.Iskovskikh and A.V.Pukhlikov, Birational
automorphisms of multi-dimensional algebraic varieties, J. Math.
Sci. {\bf 82} (1996), no. 4, 3528-3613. \vspace{0.5cm}

\noindent [K] J. Koll{\'a}r, et al., Flips and Abundance for
Algebraic Threefolds, Asterisque 211, 1993. \vspace{0.5cm}

\noindent [Kw] Y. Kawamata, A generalization of
Kodaira-Ramanujam's vanishing theorem, Math. Ann. {\bf 261}
(1982), 43-46.\vspace{0.5cm}

\noindent [P1] A.V. Pukhlikov, Birational automorphisms of
four-dimensional quintics, Invent. Math. {\bf 87} (1987), 303-329.
\vspace{0.5cm}

\noindent [P2] A.V. Pukhlikov, A note on the theorem of
V.A.Iskovskikh and Yu.I.Manin on the three-dimensional quartic,
Proc. of Steklov Math. Institute, {\bf 208} (1995), 244-254
\vspace{0.5cm}

\noindent [P3] A.V. Pukhlikov, Birational automorphisms of Fano
hypersurfaces, Invent. Math. {\bf 134} (1998), no. 2, 401-426.
\vspace{0.5cm}

\noindent [P4] A.V. Pukhlikov, Birational automorphisms of
algebraic three-folds with a pencil of del Pezzo surfaces,
Izvestiya: Mathematics {\bf 62} (1998), 115-154.
\vspace{0.5cm}

\noindent [P5] A.V. Pukhlikov, Essentials of the method of maximal
singularities, In: ``Explicit Birational Geometry of Threefolds'',
London Mathematical Society Lecture Notes Series, {\bf 281},
Cambridge University Press 2000, 73-100. \vspace{0.5cm}

\noindent [P6] A.V. Pukhlikov, Birationally rigid Fano fibrations,
Izvestiya: Mathematics {\bf 64} (2000), 131-150.
\vspace{0.5cm}

\noindent [P7] A.V. Pukhlikov, Birationally rigid Fano complete
intersections, Crelle J. f\" ur die reine und angew. Math. {\bf
541} (2001), 55-79. \vspace{0.5cm}

\noindent [P8] A.V. Pukhlikov, Birationally rigid Fano
hypersurfaces with isolated singularities, Sbornik: Mathematics
{\bf 193} (2002), no. 3; arXiv:math.AG/0106110. \vspace{0.5cm}

\noindent [Pk] J. Park, Birational maps of del Pezzo fibrations,
Crelle J. f\" ur die reine und angew. Math. {\bf 538} (2001), 213-221.
\vspace{0.5cm}

\noindent [S] V.G. Sarkisov, On conic bundle structures, Math.
USSR - Izvestiya {\bf 20} (1982), 354-390.
\vspace{0.5cm}

\noindent [Sh] V.V. Shokurov, 3-fold log flips, Izvestiya:
Mathematics {\bf 40} (1993), 93-202. \vspace{0.5cm}

\noindent [V] E. Viehweg, Vanishing theorems, Crelle J. f\" ur die reine
und angew. Math. {\bf 335} (1982), 1-8.
\vspace{0.7cm}

\begin{tabular}{ll}
Steklov Mathematics Institute & Mathematisches Institut \\
Gubkina 8 & Universit\"at Bayreuth \\
117966 Moscow & 95440 Bayreuth \\
RUSSIA & GERMANY \\
{\it e-mail}: pukh@mi.ras.ru & {\it e-mail}: pukh@btm8x5.mat.uni-bayreuth.de \\
\end{tabular}

\end{document}